\documentclass[11pt]{article}
\usepackage{amsmath,amsthm,amsfonts,latexsym,amsopn,verbatim,amscd,amssymb}
\usepackage{adjustbox}
\usepackage{lipsum}
\usepackage{cleveref}
\usepackage{enumerate}
\usepackage[left=2cm, right=2cm, top=2cm]{geometry}

\newtheorem{theorem}{Theorem}[section]

\newtheorem{Example}{Example}[section]

\newtheorem{lemma}{Lemma}[section]

\usepackage{tabularx}
\usepackage{amsmath,amsfonts,amssymb}
\usepackage{graphics, graphicx}
\usepackage{enumitem}
\usepackage{subcaption}
\usepackage{multirow}
\usepackage{authblk}
\usepackage{enumitem,amssymb}
\usepackage{graphicx,float}
\newlist{todolist}{itemize}{2}
\setlist[todolist]{label=$\square$}
\begin{document}
	\title{A hybrid approach for  singularly perturbed parabolic problem with discontinuous data}
	\date{}
\author{Nirmali Roy$^{1*}$, Anuradha Jha$^{2}$}
\date{%
    $^{1,2}$Indian Institute of Information Technology Guwahati, Bongora, India,781015\\
   *Corresponding author. E-mail: {nirmali@iiitg.ac.in}
}
	\maketitle
	\begin{sloppypar}	
\abstract{In this article, we study a two-dimensional singularly perturbed parabolic equation of the convection-diffusion type, characterized by discontinuities in the source term and convection coefficient at a specific point in the domain. These discontinuities lead to the development of interior layers. To address these layers and ensure uniform convergence, we propose a hybrid monotone difference scheme that combines the central difference and midpoint upwind schemes for spatial discretization, applied on a piecewise-uniform Shishkin mesh. For temporal discretization, we employ the Crank-Nicolson method on a uniform mesh. The resulting scheme is proven to be uniformly convergent, order achieving almost two in space and two in time. Numerical experiments validate the theoretical error estimates, demonstrating superior accuracy and convergence when compared to existing methods.}

\textbf{keywords: }{Interior layers, Hybrid method, Singularly Perturbed, Parabolic problem,  Shishkin mesh, Crank-Nicolson method}
\maketitle
\section{Introduction}\label{sec1}
We consider the following singularly perturbed parabolic differential equation with discontinuous convection coefficient and source term:
\begin{equation}
\begin{aligned}
\label{oneparaparabolic}
&\mathcal{L} y(x,t)\equiv(\epsilon y_{xx}+ ay_x-by-y_t)(x,t)=f(x,t),~(x,t)\in(\Omega^{-}\cup\Omega^{+}),\\
&y(0,t)
=p(t),~ (0,t)\in\Gamma_l =\{(0,t)|0\le t\le T\};\\
&y(1,t)=r(t),(1,t)\in\Gamma_r= \{(1,t)|0\le t\le T\};\\
&y(x,0)=q(x), (x,0)\in\Gamma_b=\{(x,0)|0\le x\le1\} ,
\end{aligned}
\end{equation}
where \(0<\epsilon\ll1\) is a singular perturbation parameter. Let \(d \in\Gamma=(0,1)\), \(\Gamma^{-}=(0,d)\), \(\Gamma^{+}=(d,1)\), \(\Gamma_{c}=\Gamma_{l}\cup\Gamma_{b}\cup\Gamma_{r}\), and \(\Omega=(0,1)\times(0,T]\), \(\Omega^{-}=(0,d)\times(0,T]\), \(\Omega^{+}=(d,1)\times(0,T]\). The convection coefficient \(a(x,t)\) and source term \(f(x,t)\) experience discontinuities at \((d,t)\in \Omega\) for all \(t\).  These functions are sufficiently smooth on \((\Omega^- \cup \Omega^+)\). We assume the jumps of \(a(x,t)\) and \(f(x,t)\) at \((d,t)\) satisfy \(|[a](d,t)|<C\) and \(|[f](d,t)|<C\), where the jump of \(\omega\) at \((d,t)\) is \([\omega](d,t)=\omega(d+,t)-\omega(d-,t)\). Moreover, \(a(x,t)\le-\alpha_1<0\) for \((x,t)\in\Omega^-\) and \(a(x,t)\ge\alpha_2>0\) for \((x,t)\in\Omega^+\), where \(\alpha_1, \alpha_2\) are positive constants. Additionally, the coefficient \(b(x, t)\) is assumed to be a sufficiently smooth function on \(\Omega\) such that \(b(x,t)\ge\beta\ge0\). 
The solution $y(x,t)$ satisfies the following interface conditions
$$[y]=0,~~~\bigg[\frac{\delta y}{\delta x}\bigg]=0, ~~\text{at} ~x=d. $$ 

Additionally, the functions \(p\), \(q\), and \(r\) are assumed to be sufficiently smooth over the domain \(\bar{\Omega}\) and must satisfy the following compatibility conditions at the two corner points \((0,0)\) and \((1,0)\):
\[q(0) = p(0),\]
\[q(1) = r(0),\]
as well as
\[\epsilon\frac{\partial^2 q(0)}{\partial x^2} + a(0)\frac{\partial q(0)}{\partial x} - b(0)q(0) - f(0,0) = \frac{\partial p(0)}{\partial t},\]
\[\epsilon\frac{\partial^2 q(1)}{\partial x^2} + a(1)\frac{\partial q(1)}{\partial x} - b(1)q(1) - f(1,0) = \frac{\partial r(0)}{\partial t}.\]
Similar compatibility conditions apply at the transition point \((d,0)\). Under these assumptions, the problem (\eqref{oneparaparabolic}) has a unique continuous solution within the domain \(\bar{\Omega}\). The solution contains strong interior layers at the points of discontinuity \((d,t)\) for all \(t \in (0,T]\) due to the discontinuity in the convection coefficient and source term.
Numerous researchers have explored singularly perturbed problems with smooth data, as highlighted in studies such as \cite{GO,KYK,PS, GKD, KS, MD, M, ER1, ZAM} as well as in case of non smooth data \cite{KM,MN,SRN06,TP,MC1}. Also in \cite{CPSH}, M. Chandru et. al. proposed a hybrid scheme to get almost second order convergence for two-parameter singularly perturbed problem with discontinuous data. These investigations have resulted in the development of various numerical methods aimed at addressing the complexities posed by such problems.   Methods like fitted meshes, finite difference schemes, and other specialized techniques have been introduced to ensure accurate and efficient solutions in these contexts.
However, when the problem data-such as source terms or convection coefficients—are discontinuous, the difficulty level increases substantially. Singularly perturbed problems with discontinuous data present unique challenges because the solution typically forms strong internal layers at point of discontinuity, which are difficult to capture numerically without a tailored approach. This growing area of research focuses on extending existing methods or developing new techniques to handle the additional complexities introduced by such discontinuities, ensuring that solutions remain robust, uniformly accurate, and computationally efficient across all parameter values.Some advances made in this area for one parameter and two parameter is in \cite{KM,MN,SRN06,TP,MC1}. Also in \cite{CPSH}, M. Chandru et. al. proposed a hybrid scheme to get almost second order convergence for two-parameter singularly perturbed problem with discontinuous data. The development of effective numerical methods for these types of problems is crucial, as they arise in many real-world applications, including chemistry \cite{JB}, lubrication theory \cite{RC}, chemical reactor theory \cite{JC}, DC motor analysis \cite{AB}.

Stynes and Roos \cite{RS} introduced a hybrid technique for steady-state singular perturbation problems (SPPs) with continuous data which was almost of order two on a Shishkin mesh. The hybrid numerical scheme for solving singularly perturbed boundary value problems (BVPs) with discontinuous convection coefficients introduced by Cen in \cite{Cen} achieves nearly second-order accuracy across the entire domain \([0,1]\) when the perturbation parameter \(\epsilon\) satisfies \(\epsilon \le N^{-1}\). In \cite{MS}, K. Mukherjeee et. al. proposed a hybrid difference scheme on Shishkin mesh in space and Euler-backward difference in time for singularly perturbed one parabolic problem with non smooth data. They achieves the convergence of order almost two in space and first order in time.

In this article, we have used a hybrid difference method to solve a singularly perturbed parabolic problem with non-smooth data, whose solution displays boundary layers and interior layers. Also, we have used the Crank-Nicolson scheme \cite{CN1} on time on a uniform mesh. The Crank-Nicolson method in time gives second order convergence in time in supremum norm. Since we aim at achieving high order of convergence, we investigate the performance of this hybrid scheme in space which combines central difference scheme, midpoint difference scheme and upwind finite difference scheme on layer adapted Shishkin mesh.

The rest of the paper is arranged as follows. In section 2, we present a \textit{priori} bounds on the analytical solution and its derivatives. The discretization of the of the problem is constructed in section 3. In section 4, truncation error of the scheme and uniform error bound is given. In section 5, we conduct few numerical experiments to validate the theoretical findings and showcase the efficiency and accuracy of the proposed scheme. The primary findings are summarized in section 6.\\
{\bf{Notation}}:  The norm used is the maximum norm given by
	$$\|y\|_{\Omega}=\max_{(x,t)\in \Omega}|y(x,t)|.$$
In the analysis we shall frequently assume that $\epsilon\le C N^{-1}$, which is the actual interest from the practical point of view. Also $C$, $C_{1}$, or $C_{2}$ are used as a positive constants that are independent of the perturbation parameter $\epsilon$, and mesh size.

 \section{Analytical aspects of the solution}
This section presents the analytical characteristics of the solution to the problem (\eqref{oneparaparabolic}). It includes a detailed examination of the solution, addressing aspects such as the existence of a solution of a singularly perturbed problem, the maximum principle, stability estimates, and a priori bounds on the derivatives of solution of problem (\ref{oneparaparabolic}). 
 \begin{lemma}
Suppose that a function $y(x,t)\in C^{0}(\bar{\Omega})\cap C^{2}(\Omega^{-}\cup \Omega^{+})$ satisfies $y(x,t)\le0,~(x,t)\in\Gamma_{c}, ~\bigg[\frac{\delta y}{\delta x}\bigg](d,t)\ge0,t>0,~~\mathcal{L}y(x,t)\ge0,~(x,t)\in \Omega^{-}\cup\Omega^{+},$ then $y(x,t)\ge0,~\forall(x,t)\in\bar{\Omega}.$
 \end{lemma}
 \begin{proof}
     The proof is given in \cite{ES}.
 \end{proof}
 A direct consequence of the previously stated maximum principle is the following stability outcome, which guarantees the boundness of the solution to the problem (\ref{oneparaparabolic}).
\begin{lemma}
    The bounds on solution $y(x,t)$ of problem \ref{oneparaparabolic} is given by 
	$$\|y\|_{\bar{\Omega}}\le \|y\|_{\Gamma_{c}}+\frac{\|f\|_{\Omega^{-}\cup\Omega^{+}}}{\theta}.$$
 where $\theta=\min\{\alpha_{1}/d,\alpha_{2}/(1-d)\}.$
\end{lemma}
\begin{proof}
     We refer to \cite{ES} for proof.
 \end{proof}
The following bounds are parameter-explicit and satisfied by the derivatives of the solution to this problem (\eqref{oneparaparabolic}). We state a results when the coefficients in problem (\eqref{oneparaparabolic}) are continuous and this result is used later on section \ref{sec4}.
\begin{lemma}
    Assume in (\ref{oneparaparabolic}) that the data $a, b\in C^{2}(\bar{\Omega}), f\in C^{2+\lambda}(\bar{\Omega})$, $\lambda\in (0,1)$ and the convection co-effcient $a(x)\ge\alpha>0, x\in \bar{\Omega}$. Also, let the initial-boundary data p , r and q be identically zero so that the compatibility conditions
    $$\frac{\delta^{i+j} f}{\delta x^{i}\delta t^{j}}(0,0)=\frac{\delta^{i+j} f}{\delta x^{i}\delta t^{j}}(1,0)=0,~~\text{for}~0\le i+2j\le2,$$
    be fulfilled. Then the solution $y(x,t)$ of the (\ref{oneparaparabolic}) lies in $C^{4+\lambda}(\bar{\Omega})$ and satisfies the estimate 
    $$\bigg|\frac{\delta^{i+j} y}{\delta x^{i}\delta t^{j}}\bigg|_{\Omega}\le C\epsilon^{-i},~~\text{for}~~0\le i+2j\le4.$$
\end{lemma}
\begin{proof}
The proof is given in \cite{NOS,CGL}.  
\end{proof}
To obtain sharper bounds on the solution, the solution $y(x,t)$ of problem (\eqref{oneparaparabolic}) is decomposed into regular and layers components as $\displaystyle y(x,t)=v(x,t)+w(x,t)$.
The regular component $v(x,t) $ satisfies the following equation:
\begin{equation}
\begin{aligned}
\label{v}
&\mathcal{L}v(x,t)=f(x,t),~~(x,t)\in\Omega^{-}\cup\Omega^{+},\\
&v(x,0)=y(x,0), ~\forall~x\in \Gamma^{-}\cup \Gamma^{+},\\
&v(0,t)=y(0,t),\\
&v(1,t)=y(1,t),t\in(0,T],\\
&v(x,0)=y(x,0),x\in\bar{\Gamma}=[0,1],
\end{aligned}
\end{equation}
and $v(d^{-},t),v(d^{+},t)$ are chosen suitably.
Furthure the regular component $v(x,t)$ can decomposed as
$$v(x,t)=\left\{
\begin{array}{ll}
\displaystyle v^{-}(x,t), & \hbox{ $(x,t)\in \Omega^-,$ }\\
v^{+}(x,t), & \hbox{ $(x,t)\in \Omega^+,$ }
\end{array}
\right.$$
where $v^{-}(x, t)$ and $v^{+}(x, t)$ are the left and right regular components respectively.

The layer component $w(x,t)$ is the solution of
\begin{equation}
\begin{aligned}
\label{lw}
&\mathcal{L}w(x,t)=0,~~(x,t)\in \Omega^{-}\cup\Omega^{+}, \\
&w(0,t)=w(1,t)=0,~\forall~ t\in(0,T],\\
&w(x,0)=0,~\forall~x\in \bar{\Omega},\\
&[w](d,t)=-[v](d,t),~\bigg[\frac{\delta w}{\delta x}\bigg](d,t)=-\bigg[\frac{\delta v}{\delta x}\bigg](d,t),~t\in (0,T].
\end{aligned}
\end{equation}
Hence, $w(d^{-},t)=y(d^{-},t)-v(d^{-},t)$ and $w(d^{+},t)=y(d^{+},t)-v(d^{+},t).$
Further, the layer component $w(x,t)$ is decomposed as
$$w(x,t)=\left\{
\begin{array}{ll}
\displaystyle w^{-}(x,t), & \hbox{ $(x,t)\in \Omega^-,$ }\\
w^{+}(x,t), & \hbox{ $(x,t)\in \Omega^+.$ }
\end{array}
\right.
$$
Hence, the unique solution $y(x,t)$ to problem (\eqref{oneparaparabolic}) is written as
$$y(x,t)=\left\{
\begin{array}{ll}
\displaystyle (v^{-}+w^{-})(x,t),\hspace{2.5cm} (x,t)\in\Omega^{-},\\
(v^{-}+w^{-})(d^{-},t)\\=(v^{+}+w^{+})(d^{+},t), ~~\hspace{1.5cm} (x,t)=(d,t),\forall~ t\in(0,T],\\
(v^{+}+w^{+})(x,t),\hspace{2.5cm}(x,t)\in\Omega^{+}.
\end{array}
\right.$$

 \begin{lemma}\cite{LSU}
 \label{bounds1}
    The bounds on the components $v(x,t), w(x,t)$ and their derivatives are given as follows $0\le i\le 3$ and $0\le i+2j\le 4$:
    $$\bigg\|\frac{\delta^{i+j}v}{\delta x^{i}\delta t^{j}}\bigg\|_{\Omega^{-}\cup
		\Omega^{+}}\le C,~~~~~\bigg\|\frac{\delta^{4}v}{\delta x^{4}}\bigg\|_{\Omega^{-}\cup
		\Omega^{+}}\le C \epsilon^{-1},$$
  
     $$\bigg\|\frac{\delta^{i+j}w}{\delta x^{i}\delta t^{j}}\bigg\|_{\Omega^{-}\cup
		\Omega^{+}} \le C \epsilon^{-i}\left\{
	\begin{array}{ll}
	\displaystyle e^{-(d-x)\alpha_{1}/\epsilon}, & \hbox{ $(x, t)\in \Omega^{-},$ }\\
	e^{-(x-d)\alpha_{2}/\epsilon}, & \hbox{ $(x, t)\in \Omega^{+}$,}
	\end{array}
	\right.$$
	$$\bigg\|\frac{\delta^{4}w}{\delta x^{4}}\bigg\|_{\Omega^{-}\cup
		\Omega^{+}}\le C \epsilon^{-4} \left\{
	\begin{array}{ll}
	\displaystyle e^{-(d-x)\alpha_{1}/\epsilon}, & \hbox{ $(x, t)\in \Omega^{-},$ }\\
	e^{-(x-d)\alpha_{2}/\epsilon}, & \hbox{ $(x, t)\in \Omega^{+}$,}
	\end{array}
	\right.$$
and $C$ is a constant independent of $\epsilon$.
 \end{lemma}

 \section{Semi discretization}
 In this section, we begin by demonstrating the semi-discretization of the problem in the temporal direction on a uniform mesh using the Crank-Nicolson method \cite{CN1}. The Crank-Nicolson method is second-order accurate in time, meaning it provides a higher level of precision compared to first-order methods like backward Euler.
 For a fixed time $T$, the interval $[0,T]$ is partitioned uniformly as $\Lambda^{M}=\{t_{j}=j\Delta t:j=0,1,\ldots,M, \Delta t=\frac{T}{M}\}$. The semi-discretization yields the following system of linear ordinary differential equations:
\begin{align*}
&\epsilon Y^{j+\frac{1}{2}}_{xx}(x)+ a^{j+\frac{1}{2}}(x)Y^{j+\frac{1}{2}}_{x}(x)-b^{j+\frac{1}{2}}(x)Y^{j+\frac{1}{2}}(x)=f^{j+\frac{1}{2}}(x)+\frac{Y^{j+1}(x)-Y^{j}(x)}{\Delta t},\\
&\hspace{3.3in} x\in  \Omega^{-}\cup\Omega^{+},~ 0\le j\le M-1,\\
&Y^{j+1}(0)=y(0,t_{j+1}),~~~Y^{j+1}(1)=y(1,t_{j+1}),~~0\le j\le M-1,\\
&Y^{0}(x)=y(x,0),~~x\in\Gamma_{b},
\end{align*}
where $Y^{j+1}(x)$ is the approximation of $y(x,t_{j+1})$ of problem (\eqref{oneparaparabolic}) at $(j+1)$-th time level and $\displaystyle Y^{j+\frac{1}{2}}=\frac{Y^{j+1}(x)+Y^{j}(x)}{2} $.
Upon simplification, we have
\begin{equation}
\label{semi-dis}
\left.
\begin{array}{ll}
\displaystyle \mathcal{\tilde{L}}Y^{j+1}(x)=g(x,t_{j}), & \hbox{ $x\in  \Gamma^{-}\cup\Gamma^{+},~ 0\le j\le M-1$, }\\
Y^{j+1}(0)=y(0,t_{j+1}), & \hbox{ $0\le j\le M-1$,}\\Y^{j+1}(1)=y(1,t_{j+1}), & \hbox{ $0\le j\le M-1$,}\\
Y^{0}(x)=y(x,0), & \hbox{ $x\in\Gamma,$ }
\end{array}
\right\}
\end{equation}
where the operator $\mathcal{\tilde{L}}$ is defined as\\
$$\mathcal{\tilde{L}}\equiv \epsilon \frac{d^{2}}{d x^{2}}+ a^{j+\frac{1}{2}}\frac{d}{d x}-c^{j+\frac{1}{2}}I,$$\\
and
\begin{align*}
    g(x,t_{j})&=2f^{j+\frac{1}{2}}(x)-\epsilon Y^{j}_{xx}(x)-a^{j+\frac{1}{2}}(x)Y_{x}^{j}(x)+d^{j+\frac{1}{2}}(x)Y^{j}(x),\\
c^{j+\frac{1}{2}}(x)&=b^{j+\frac{1}{2}}(x)+\frac{2}{\Delta t},\\
d^{j+\frac{1}{2}}(x)&=b^{j+\frac{1}{2}}(x)-\frac{2}{\Delta t}.
\end{align*}
The error in temporal semi-discretization is defined by 
$e^{j+1}= y(x,t_{j+1})- \hat{Y}^{j+1}(x)$, where $y(x,t_{j+1})$ is the solution of  problem (\eqref{oneparaparabolic}). $ \hat{Y}^{j+1}(x)$  is the solution of semi-discrete scheme \eqref{semi-dis}, when $y(x,t_{j})$ is taken instead of $Y^{j}$ to find solution at $(x,t_{j+1})$.
 \begin{theorem}
     The local truncation error $T_{j+1}=\mathcal{\tilde{L}}(e^{j+1})$ satisfies 
	$$\|T_{j+1}\|\le C(\Delta t)^{3},~~~0\le j\le M-1.$$
 \end{theorem}
 \begin{proof}
The result can be proved following the approach given in \cite{KTK}.
\end{proof}
\begin{theorem}
	The global  error $E_{\epsilon}^{j+1}=y(x,t_{j+1})-Y^{j+1}(x)$ is estimated as 
	$$\|E_{\epsilon}^{j+1}\|\le C(\Delta t)^{2},~~~0\le j\le M-1,$$
	where $Y^{j+1}(x)$ is the solution of \eqref{semi-dis}.
\end{theorem}
\begin{proof}
	For proof, we follow the approch in \cite{KTK}.
\end{proof}
\subsection{Discretization in space}
We understand from the relevant literature on singular perturbation that a uniform mesh on standard difference methods is insufficient for achieving uniform convergence for solutions of the problem (\eqref{oneparaparabolic}) because of the existence of layer regions. The differential equation (\eqref{oneparaparabolic}) exhibits strong interior layers. To resolve these layers, we have discretized this problem by using a piece-wise uniform Shishkin mesh.\\
Let $N\ge8$ be a multiple of 4. The interior points of the spatial mesh be denoted by $\Gamma^{N}=\{x_{i}:1\le i \le \frac{N}{2}-1\}\cup\{x_{i}:\frac{N}{2}+1\le i \le N-1\}.$
The $\bar{\Gamma}^{N}=\{x_{i}\}_{0}^{N}$ denote the mesh points with $x_{0}=0,x_{N}=1$ and the point of discontinuity at point $x_{\frac{N}{2}}=d.$ We also introduce the notation $\Gamma^{N-}=\{x_{i}\}_{0}^{\frac{N}{2}-1}, \Gamma^{N+}=\{x_{i}\}^{N-1}_{\frac{N}{2}+1}$, $\Omega^{N-}=\Gamma^{N-}\times\Lambda^{M}$, $\Omega^{N+}=\Gamma^{N+}\times\Lambda^{M}$, $\Omega^{N,M}=\Gamma^{N}\times\Lambda^{M}, \bar{\Omega}^{N,M}=\bar{\Gamma}^{N}\times\Lambda^{M}$, $\Gamma_{c}^{N}=\bar{\Omega}^{N,M}\cap\Gamma_{c}$, $\Gamma^{N}_{l}=\Gamma_{c}^{N}\cap\Gamma_{l}$, $\Gamma_{b}^{N}=\Gamma_{c}^{N}\cap\Gamma_{b}$ and $\Gamma^{N}_{r}=\Gamma_{c}^{N}\cap\Gamma_{r}$.
The domain $[0,1]$ is subdivided into four sub-intervals as
\[\bar{\Gamma}=[0,d-\tau_{1]}]\cup[d-\tau_{1},d]\cup[d,d+\tau_{2]}]\cup[d+\tau_{2},1]= \bigcup_{i=1}^{4}I_{i}.\]
The transition point in $\Gamma^{N}$ is:
$$\tau_{1}=\min\bigg\{\frac{d}{2},\frac{2\epsilon}{\alpha}\ln N\bigg\}.$$
$$\tau_{2}=\min\bigg\{\frac{1-d}{2},\frac{2\epsilon}{\alpha}\ln N\bigg\}.$$
The mesh points are given by
\begin{equation*}
x_{i}=\left\{
\begin{array}{ll}
\displaystyle \frac{4i(d-\tau_{1})}{N}, & 0\le i \le \frac{N}{4}, \vspace{.2cm} \\
\displaystyle (d-\tau_{1})+\frac{4\tau_{1}(i-\frac{N}{4})}{N}, & \frac{N}{4}\le i \le \frac{N}{2}, \vspace{.2cm}\\
\displaystyle d+\frac{4(i-\frac{N}{2})\tau_{2}}{N}, & \frac{N}{2}\le i \le \frac{3N}{4},\vspace{.2cm} \\

\displaystyle (d+\tau_{2})+\frac{4(i-\frac{3N}{4})(1-d-\tau_{2})}{N}, &  \frac{3N}{4}\le i \le N.  
\end{array}
\right.
\end{equation*}
Let $h_{i}=x_{i}-x_{i-1}$ denote the step size and $H_{1}=\frac{4(d-\tau_{1})}{N}, H_{2}=\frac{4\tau_{1}}{N},H_{3}=\frac{4\tau_{2}}{N},H_{4}=\frac{4(1-d-\tau_{2})}{N}$ denotes the step sizes in the interval $I_{1},I_{2},I_{3},I_{4}$, respectively. Also, we use $~h=H_{2}=H_{3}$, whenever $\tau_{1}=\tau_{2}.$ Also $\hbar=\frac{h_{i}+h_{i+1}}{2}.$ \\
We observe that $$N^{-1}\le \{H_{1},H_{4}\}\le 4N^{-1}, h+H_{1}=h+H_{4}\le 4N^{-1}.$$
We have used a hybrid difference approach for the discretization of the problem (\eqref{oneparaparabolic}) in the spatial variable. On the above piece-wise uniform Shishkin mesh, this hybrid difference scheme comprises of a central, and a midpoint scheme :
 \small
 \[\mathcal{L}_{h}^{N,M}Y^{j+1}(x_{i})=\left\{
\begin{array}{ll}
\displaystyle \mathcal{L}_{m}^{N,M}Y^{j+1}(x_{i}), & x_{i}\in(0,d-\tau_{1}]\cup[d+\tau_{2},1),\\
\displaystyle \mathcal{L}_{c}^{N,M}Y^{j+1}(x_{i}), & x_{i}\in(d-\tau_{1},d)\cup(d,d+\tau_{2}),\\
\mathcal{L}_{t}^{N,M}Y^{j+1}(x_{\frac{N}{2}}), &  x_{i}=d.
\end{array}
\right.\]

The fully discretize scheme is given by
\begin{equation}
\label{discrete problem}
    \mathcal{L}_{h}^{N,M}Y^{j+1}(x_{i})\equiv r_{i}^{-}Y^{j+1}(x_{i-1})+r_{i}^{c}Y^{j+1}(x_{i})+r_{i}^{+}Y^{j+1}(x_{i+1})=g^{j+1}(x_{i}).
\end{equation}
At the point of discontinuity $x_{N/2}=d$, we have used five point second order difference scheme
\begin{equation}
    \begin{aligned}
    \label{five point scheme}
        \mathcal{L}^{N,M}_{t}Y^{j+1}(x_{\frac{N}{2}})=\frac{-Y^{j+1}(x_{N/2+2})+4Y^{j+1}(x_{\frac{N}{2}+1})-3Y^{j+1}(x_{\frac{N}{2}})}{2H_{4}}\\-\frac{Y^{j+1}(x_{\frac{N}{2}-2})-4Y^{j+1}(x_{\frac{N}{2}-1})+3Y^{j+1}(x_{\frac{N}{2}})}{2H_{3}}=0.
    \end{aligned}
\end{equation}
Associated with each of these finite difference operators, we have the following finite difference scheme
\begin{equation*}
\begin{aligned}
\label{schemes}
 \mathcal{L}_{m}^{N,M}Y^{j+1}(x_{i})\equiv \epsilon \delta^{2}Y^{j+1}(x_{i})+ \Bar{a}^{j+\frac{1}{2}}(x_{i})D^{*}Y^{j+1}(x_{i})-\Bar{b}^{j+\frac{1}{2}}(x_{i})\bar{Y}^{j+1}(x_{i}),\\
  \mathcal{L}_{c}^{N,M}Y^{j+1}(x_{i})\equiv \epsilon \delta^{2}Y^{j+1}(x_{i})+ a^{j+\frac{1}{2}}(x_{i})D^{0}Y^{j+1}(x_{i})-b^{j+\frac{1}{2}}(x_{i})Y^{j+1}(x_{i}),
 \end{aligned}
\end{equation*}
and
\[g^{j+1}(x_{i})=\left\{
\begin{array}{ll}
\displaystyle 2\bar{f}^{j+\frac{1}{2}}(x_{i})-\epsilon \delta_{x}^{2}Y^{j}(x_{i})- \bar{a}^{j+\frac{1}{2}}(x_{i})D_{x}^{*}Y^{j}(x_{i})+\bar{d}^{j+\frac{1}{2}}(x)Y^{j}(x_{i}), & \mathcal{L}_{h}^{N,M}=\mathcal{L}_{m}^{N,M}, \vspace{0.2cm}\\
\displaystyle 2f^{j+\frac{1}{2}}(x_{i})-\epsilon \delta_{x}^{2}Y^{j}(x_{i})- a^{j+\frac{1}{2}}(x_{i})D_{x}^{0}Y^{j}(x_{i})+d^{j+\frac{1}{2}}(x)Y^{j}(x_{i}), & \mathcal{L}_{h}^{N,M}=\mathcal{L}_{c}^{N,M}, \vspace{0.2cm}
\end{array}
\right.\]
 where,
 \[D^{+}Y^{j+1}(x_{i})=\frac{Y^{j+1}(x_{i+1})-Y^{j+1}(x_{i})}{x_{i+1}-x_{i}},~~~~ D^{-}Y^{j+1}(x_{i})=\frac{Y^{j+1}(x_{i})-Y^{j+1}(x_{i-1})}{x_{i}-x_{i-1}},\]
 \[D^{0}Y^{j+1}(x_{i})=\frac{Y^{j+1}(x_{i+1})-Y^{j+1}(x_{i-1})}{x_{i+1}-x_{i-1}},~~~~\delta^{2}Y^{j+1}(x_{i})=\frac{2(D^{+}Y^{j+1}(x_{i})-D^{-}Y^{j+1}(x_{i}))}{x_{i+1}-x_{i-1}},\]
\[D^{*}Y^{j+1}(x_{i})=\left\{
\begin{array}{ll}
\displaystyle D^{-}Y^{j+1}(x_{i}), & i<\frac{N}{2}, \vspace{0.2cm}\\
D^{+}Y^{j+1}(x_{i}), &  i>\frac{N}{2},
\end{array}
\right.\] and
$\displaystyle \bar{w}^{j+1}(x_{i})=\frac{w^{j+1}(x_{i})+w^{j+1}(x_{i-1})}{2}$ in $\Omega^{N-}$,$~~\displaystyle \bar{w}^{j+1}(x_{i})=\frac{w^{j+1}(x_{i})+w^{j+1}(x_{i+1})}{2}$ in $\Omega^{N+}$.\\
As the matrix associated with the above equations \eqref{discrete problem} is not an M-matrix. We transform the equation (\ref{five point scheme}) to establish the monotonicity property by estimating $Y_{\frac{N}{2}-2}^{j+1}$ and $Y_{\frac{N}{2}+2}^{j+1}$ for $\mathcal{L}_{t}^{N,M}Y^{j+1}(x_{i})$ from the operator $\mathcal{L}_{c}^{N,M}$, we get\\ 
$Y_{N/2-2}^{j+1}=\frac{2h^2}{2\epsilon- h a_{N/2-1}^{j+\frac{1}{2}}}\bigg\{ g_{N/2-1}^{j+1}-Y_{N/2-1}^{j+1}\bigg(\frac{-2\epsilon-h^2(b_{N/2-1}^{j+\frac{1}{2}}+\frac{2}{\Delta t})}{h^2}\bigg)-Y_{N/2}^{j+1}\bigg(\frac{2\epsilon+ ha_{N/2-1}^{j+\frac{1}{2}}}{2h^2}\bigg)\bigg\}$\\
$Y_{N/2+2}^{j+1}=\frac{2h^2}{2\epsilon+ h a_{N/2+1}^{j+\frac{1}{2}}}\bigg\{ g_{N/2+1}^{j+1}-Y_{N/2+1}^{j+1}\bigg(\frac{-2\epsilon-h^2(b_{N/2+1}^{j+\frac{1}{2}}+\frac{2}{\Delta t})}{h^2}\bigg)-Y_{N/2}^{j+1}\bigg(\frac{2\epsilon- ha_{N/2+1}^{j+\frac{1}{2}}}{2h^2}\bigg)\bigg\}$\\
Inserting the above expression in $\mathcal{L}^{N,M}_{t}Y^{j+1}(x_{N/2})$, results a tridiagonal system of equations of the following form:\\
\begin{align*}
\mathcal{L}^{N,M}_{t}Y^{j+1}_{N/2}&=\bigg(\frac{2\epsilon-h a_{N/2+1}^{j+\frac{1}{2}}}{2\epsilon+h a_{N/2+1}^{j+\frac{1}{2}}}-6+\frac{2\epsilon+h a_{N/2-1}^{j+\frac{1}{2}}}{2\epsilon-ha_{N/2-1}^{j+\frac{1}{2}}}\bigg)Y^{j+1}_{N/2}+\bigg(\frac{-4\epsilon-2h^2 (b_{N/2+1}^{j+\frac{1}{2}}+\frac{2}{\Delta t})}{2\epsilon+h a_{N/2+1}^{j+\frac{1}{2}}}+4\bigg)Y^{j+1}_{N/2+1}\\&+\bigg(\frac{-4\epsilon-2h^2 (b_{N/2-1}^{j+\frac{1}{2}}+\frac{2}{\Delta t})}{2\epsilon-h a_{N/2-1}^{j+\frac{1}{2}}}+4\bigg)Y_{N/2-1}^{j+1}.\\&=\frac{2h^2 g^{j+1}_{N/2-1}}{2\epsilon-h a_{N/2-1}^{j+\frac{1}{2}}}+\frac{2h^2 g_{N/2+1}^{j+1}}{2\epsilon+h a_{N/2+1}^{j+\frac{1}{2}}}.
\end{align*}
The element of the system matrix $\mathcal{L}_{h}^{N,M}$ are as follows:\\
$\text{In}~\Omega^{-},~\text{if}~ \mathcal{L}_{h}^{N,M}=\mathcal{L}_{m}^{N,M},$
$$r_{i}^{-}=\frac{\epsilon}{h_{i}\hbar_{i}}-\frac{ \Bar{a}^{j+\frac{1}{2}}(x_{i})}{h_{i}}-\bigg(\frac{\bar{b}^{j+\frac{1}{2}}(x_{i})}{2}+\frac{1}{\Delta t}\bigg), ~r_{i}^{+}=\frac{\epsilon}{h_{i+1}\hbar_{i}},~r_{i}^{c}=-r_{i}^{-}-r_{i}^{+}-\frac{1}{2}\bigg(\bar{b}^{j+\frac{1}{2}}(x_{i})+\frac{2}{\Delta t}\bigg).$$
$\text{In}~\Omega^{+},~\text{if}~ \mathcal{L}_{h}^{N,M}=\mathcal{L}_{m}^{N,M},$
$$\text{In}~\Omega^{+},~r_{i}^{-}=\frac{\epsilon}{h_{i}\hbar_{i}}, ~r_{i}^{+}=\frac{\epsilon}{h_{i+1}\hbar_{i}}-\frac{ \Bar{a}^{j+\frac{1}{2}}(x_{i})}{h_{i+1}}-\bigg(\frac{\bar{b}^{j+\frac{1}{2}}(x_{i})}{2}+\frac{1}{\Delta t}\bigg),~r_{i}^{c}=-r_{i}^{-}-r_{i}^{+}-\frac{1}{2}\bigg(\bar{b}^{j+\frac{1}{2}}(x_{i})+\frac{2}{\Delta t}\bigg).$$
$\text{In}~\Omega,~\text{if}~ \mathcal{L}_{h}^{N,M}=\mathcal{L}_{c}^{N,M},$
$$r_{i}^{-}=\frac{\epsilon}{h_{i}\hbar_{i}}-\frac{ a^{j+\frac{1}{2}}(x_{i})}{2\hbar_{i}}, ~r_{i}^{+}=\frac{\epsilon}{h_{i+1}\hbar_{i}}+\frac{ a^{j+\frac{1}{2}}(x_{i})}{2\hbar_{i}},~r_{i}^{c}=-r_{i}^{-}-r_{i}^{+}-\bigg(b^{j+\frac{1}{2}}(x_{i})+\frac{2}{\Delta t}\bigg).$$

$$\text{At $x_{N/2}=d$},~r_{N/2}^{-}=\bigg(\frac{-4\epsilon-2h^2 (b_{N/2-1}^{j+\frac{1}{2}}+\frac{2}{\Delta t})}{2\epsilon-h a_{N/2-1}^{j+\frac{1}{2}}}+4\bigg), ~r_{N/2}^{+}=\bigg(\frac{-4\epsilon-2h^2 (b_{N/2+1}^{j+\frac{1}{2}}+\frac{2}{\Delta t})}{2\epsilon+h a_{N/2+1}^{j+\frac{1}{2}}}+4\bigg),$$ 
$$~r_{N/2}^{c}=\bigg(\frac{2\epsilon-h a_{N/2+1}^{j+\frac{1}{2}}}{2\epsilon+h a_{N/2+1}^{j+\frac{1}{2}}}-6+\frac{2\epsilon+h a_{N/2-1}^{j+\frac{1}{2}}}{2\epsilon-ha_{N/2-1}^{j+\frac{1}{2}}}\bigg).$$
To find the error estimates for the scheme \eqref{discrete problem} defined  above, we first decompose the discrete solution $Y(x_{i},t_{j+1})$ into the discrete regular and discrete singular components. \\
Let
\[Y^{j+1}(x_{i})=V^{j+1}(x_{i})+W^{j+1}(x_{i}).\]
The discrete regular components further decomposed as 
$$V^{j+1}(x_{i})=\left\{
\begin{array}{ll}
\displaystyle V^{-(j+1)}(x_{i}), & \hbox{ $(x_{i},t_{j+1})\in \Omega^{N-},$ }\\
V^{+(j+1)}(x_{i}),& \hbox{ $(x_{i},t_{j+1})\in \Omega^{N+},$ }
\end{array}
\right.$$
where $V^{-(j+1)}(x_{i})$ and $V^{+(j+1)}(x_{i})$ approximate $v^{-}(x_{i},t_{j+1})$ and $v^{+}(x_{i},t_{j+1})$ respectively. They satisfy the following equations:
\begin{equation}
\label{LV}
\begin{aligned}
	\left.
	\begin{array}{ll} 
\mathcal{L}^{N,M}_{h}V^{-(j+1)}(x_{i})=g(x_{i},t_{j+1}),~\forall (x_{i},t_{j+1})\in\Omega^{N-},\\V^{-(j+1)}(x_{0})=v^{-}(0,t_{j+1}), V^{-(j+1)}(x_{\frac{N}{2}})=v^{-}(d-,t_{j+1}),\\
\mathcal{L}^{N,M}_{h}V^{+(j+1)}(x_{i})=g(x_{i},t_{j+1}),~\forall (x_{i},t_{j+1})\in\Omega^{N+},\\
V^{+(j+1)}(x_{\frac{N}{2}})=v^{+}(d+,t_{j+1}), V^{+(j+1)}(x_{N})=v^{-}(1,t_{j+1}).
\end{array}
\right\}
\end{aligned}
\end{equation}
The discrete singular components $W^{j+1}(x_{i})$ is also decomposed as:
$$W^{j+1}(x_{i})=\left\{
\begin{array}{ll}
\displaystyle W^{-(j+1)}(x_{i}), & \hbox{ $(x_{i},t_{j+1})\in \Omega^{N-},$ }\\
W^{+(j+1)}(x_{i}),& \hbox{ $(x_{i},t_{j+1})\in \Omega^{N+},$ }
\end{array}
\right.$$

where $W^{-(j+1)}(x_{i})$ approximates the layer components $w^{-}(x_{i},t_{j+1})$ and $w^{+}(x_{i},t_{j+1})$ respectively.   
These components satisfy the following equations
\begin{equation}
\label{LW-}
\begin{aligned}
	\left.
	\begin{array}{ll} 
\mathcal{L}^{N,M}_{h}W^{-(j+1)}(x_{i})=0,~\forall (x_{i},t_{j+1})\in\Omega^{N-}, \\ W^{-(j+1)}(x_{0})=w^{-}(0,t_{j+1}), W^{-(j+1)}(x_{\frac{N}{2}})=w^{-}(d,t_{j+1}),
\end{array}
\right\}
\end{aligned}
\end{equation}
\begin{equation}
\label{LW+}
\begin{aligned}
	\left.
	\begin{array}{ll} 
\mathcal{L}^{N,M}_{h}W^{+(j+1)}(x_{i})=0,~\forall ~(x_{i},t_{j+1})\in\Omega^{N+},\\
W^{+(j+1)}(x_{\frac{N}{2}})=w^{+}(d,t_{j+1}), W^{+(j+1)}(x_{N})=w^{+}(1,t_{j+1}).
\end{array}
\right\}
\end{aligned}
\end{equation}
Hence, the complete discrete solution $U^{j+1}(x_{i})$ is defined by
$$Y^{j+1}(x_{i})=\left\{
\begin{array}{ll}
\displaystyle (V^{-(j+1)}+W^{-(j+1)})(x_{i}),\hspace{1.1cm} (x_{i},t_{j+1})\in\Omega^{N-},\\
(V^{-(j+1)}+W^{-(j+1)})(d-,t)\\
=(V^{+(j+1)}+W^{+(j+1)})(d+,t),\hspace{0.5cm}(x_{i},t_{j+1})=(d,t_{j+1}),\\
(V^{+(j+1)}+W^{+(j+1)}),\hspace{1.1cm}(x_{i},t_{j+1})\in\Omega^{N+}.
\end{array}
\right.$$
\begin{lemma}
Let $$N(\ln N)^{-1}>\frac{4\|a\|}{\alpha},~~ 2N\|a\|\ge\|b\|+\frac{2M}{T},$$
where $\alpha=\min\{\alpha_{1}, \alpha_{2}\}$. The operator defined by (\ref{discrete problem}) satisfies a discrete maximum principle, i.e., 
    suppose that a mesh function $Y(x_{i},t_{j})$ satisfies $Y(x_{i},t_{j})\le 0,~\forall (x_{i},t_{j})\in \Gamma_{c}^{N}$ and $\mathcal{L}_{t}^{N,M}Y(x_{N/2})\ge0 ~\forall~j=0,\ldots,M$. If $\mathcal{L}_{h}^{N,M}Y(x_{i},t_{j})\ge0$ for all $(x_{i},t_{j})\in \Omega^{N-}\cup\Omega^{N+}$ then $Y(x_{i},t_{j})\le0,~\forall  (x_{i},t_{j})\in\bar{\Omega}^{N,M}$.
\end{lemma}
\begin{proof}
  We follow a proof analogous to one given in \cite{MS}.
\end{proof}
 A direct implication of the discrete maximum principle is the resulting stability bounds.
\begin{lemma}
    If $Y(x_{i},t_{j+1}), (x_{i},t_{j+1})\in\bar{\Omega}^{N,M}$ is a mesh function  satisfying the difference scheme \eqref{discrete problem}, then $ \|Y\|_{\bar{\Omega}^{N,M}}\le C$.
    
\end{lemma}
\begin{proof}
    Define the mesh function $\psi(x_{i},t_{j+1})=-\phi(x_{i},t_{j+1})\pm Y(x_{i},t_{j+1}),$
    where $\phi(x_{i},t_{j+1})=\max \bigg\{|Y^{j+1}(0)|,|Y^{j+1}(1)|,\frac{\|g\|_{\bar{\Omega}^{N,M}}}{\theta}\bigg\},$ where $\theta=\min\{\alpha_{1}/d,\alpha_{2}/(1-d)\}.$
    Clearly $\psi(x_{i},t_{j+1})\le 0~\forall~(x_{i},t_{j+1})\in \Gamma_{c}^{N}$ and 
    $\mathcal{L}^{N,M}_{c}\psi(x_{i},t_{j+1})=-b(x_{i},t_{j+\frac{1}{2}})\phi(x_{i},t_{j+1})\pm \mathcal{L}^{N,M}_{c}Y(x_{i},t_{j+1})\ge0,$
    
    $\mathcal{L}^{N,M}_{m}\psi(x_{i},t_{j+1})=-\Bar{b}(x_{i},t_{j+\frac{1}{2}})\bar{\phi}(x_{i},t_{j+1})\pm \mathcal{L}^{N,M}_{m}Y(x_{i},t_{j+1})\ge0.$
     At the point of discontinuity, 
    $\mathcal{L}^{N,M}_{t}\psi(x_{\frac{N}{2}},t_{j+1})\ge0.$\\
    Therefore, $\psi(x_{i},t_{j+1})\le0,~~\forall(x_{i},t_{j+1})\in\bar{\Omega}^{N,M}.$
    This leads to the required result
    $$\|Y\|_{\bar{\Omega}^{N,M}}\le C$$
\end{proof}
\section{Truncation error analysis}\label{sec4}
\begin{lemma}
\label{discrete singular}
    The layer component $W^{-}(x_{i},t_{j+1})$ and $ W^{+}(x_{i},t_{j+1})$ satisfy the following bounds
	$$|W^{-}(x_{i},t_{j+1})|\le C R_{i},\quad R_{i}=\prod_{n=i+1}^{\frac{N}{2}}\bigg(1+\frac{\alpha h_{n}}{\epsilon}\bigg)^{-1},~~R_{0}=C_{1},~i=1,\ldots,\frac{N}{2}-1,$$
	$$|W^{+}(x_{i},t_{j+1})|\le C Q_{i},\quad Q_{i}=\prod_{n=\frac{N}{2}+1}^{i}\bigg(1+\frac{\alpha h_{n}}{\epsilon}\bigg)^{-1},~~Q_{\frac{N}{2}}=C_{1},~i=\frac{N}{2}+1,\ldots,N,$$
\end{lemma}
\begin{proof}
Define the barrier function for the left layer term as
\[\eta^{-(j+1)}_{i}= -C R_{i} \pm W^{-(j+1)}(x_{i}),\quad  0\leq i \leq \frac{N}{2}-1,~0\le j\le M-1 .\]	
where 
\[R_{i}^{(j+1)}=\left\{
\begin{array}{ll}
\displaystyle \prod_{n=i+1}^{\frac{N}{2}}\bigg(1+\frac{\alpha h_{n}}{\epsilon}\bigg)^{-1}, & 1\le i\le \frac{N}{2}-1,\vspace{0.2cm}\\
1, &  i=0.
\end{array}
\right.\]
 In $0\le i \le \frac{N}{4}$:
\begin{eqnarray*}
    -\mathcal{L}_{m}^{N,M} R_{i}&=&- \frac{2\alpha}{h_{i}+h_{i+1}}(R_{i}-R_{i-1})-\bar{a}^{j+\frac{1}{2}}(x_{i})R_{i-1}\frac{\alpha}{\epsilon}+\bar{d}^{j+\frac{1}{2}}(x_{i})R_{i}\\
 &\ge&-\frac{\alpha}{\epsilon}R_{i-1}\bigg(\frac{2\alpha h_{i}}{h_{i}+h_{i+1}}+\bar{a}^{j+\frac{1}{2}}(x_{i})\bigg)\\&\ge& \frac{C R_{i}}{\epsilon+\alpha h_{i}}\ge 0
\end{eqnarray*}
\\
Now, $\frac{N}{4}\le i \le \frac{N}{2},$
 \begin{eqnarray*}
 -\mathcal{L}_{c}^{N,M} R_{i}&=&- \frac{\alpha}{h}(R_{i}-R_{i-1})-a^{j+\frac{1}{2}}(x_{i})R_{i-1}\bigg(\frac{\alpha}{\epsilon}+\frac{\alpha^{2}h}{2\epsilon^2}\bigg)+d^{j+\frac{1}{2}}(x_{i})R_{i}\\&\ge& -\frac{\alpha}{\epsilon}R_{i-1}(\alpha+a^{j+\frac{1}{2}}(x_{i}))-a^{j+\frac{1}{2}}(x_{i})\frac{\alpha^{2}h}{2\epsilon^{2}}R_{i-1}\\&\ge&\frac{C R_{i}}{\epsilon+\alpha h}\ge 0
 \end{eqnarray*}
Also for large enough $C$, $\eta^{-(j+1)}_{l,0}\leq 0,\eta^{-(0)}_{i}\leq 0$ and $\eta^{-(j+1)}_{l,N/2}\leq 0$.
 So, by discrete maximum principle, we obtain  \[\eta^{-(j+1)}_{i} \leq  0 \implies |W^{-(j+1)}(x_{i})|\leq C\prod_{n=i+1}^{\frac{N}{2}}\bigg(1+\frac{\alpha h_{n}}{\epsilon}\bigg)^{-1}, \; 1\leq i \leq \frac{N}{2}-1,~0\le j\le M-1 .\]
Now, to prove the bound for $W^{+(j+1)}(x_{i})$, we define the barrier function for the right layer component  as
\[\eta^{+(j+1)}_{i}= -CQ_{i} \pm W^{+(j+1)}(x_{i}),\quad  \frac{N}{2}+1\leq i \leq N-1,~0\le j\le M-1.\]	
where 
\[Q_{i}=\left\{
\begin{array}{ll}
\displaystyle \prod_{n=\frac{N}{2}+1}^{i}\bigg(1+\frac{\alpha h_{n}}{\epsilon}\bigg)^{-1}, & \frac{N}{2}+1\le i< N-1, \vspace{0.2cm}\\
1, &  i=\frac{N}{2}.
\end{array}
\right.\]
Now, $\frac{N}{2}+1\le i\le \frac{3N}{4}-1$
\begin{eqnarray*}
	 -\mathcal{L}_{c}^{N,M} Q_{i}&=&\frac{\alpha}{h}(Q_{i+1}-Q_{i})+a^{j+\frac{1}{2}}(x_{i})Q_{i+1}\bigg(\frac{\alpha}{\epsilon}+\frac{\alpha^{2}h}{2\epsilon^2}\bigg)+d^{j+\frac{1}{2}}(x_{i})Q_{i}\\&\ge& \frac{\alpha}{\epsilon}Q_{i+1}\big(a^{j+\frac{1}{2}}(x_{i})-\alpha\big)+a^{j+\frac{1}{2}}(x_{i})\frac{\alpha^{2}h}{2\epsilon^{2}}Q_{i+1}\\&\ge&\frac{C Q_{i}}{\epsilon+\alpha h}\ge 0
\end{eqnarray*}
In $\frac{3N}{4}\le i \le N-1$:
\begin{eqnarray*}
	 -\mathcal{L}_{m}^{N,M} Q_{i}&=& \frac{2\alpha}{h_{i}+h_{i+1}}(Q_{i+1}-Q_{i})+\bar{a}^{j+\frac{1}{2}}(x_{i})Q_{i}\bigg(\frac{\alpha}{\epsilon}\bigg)+\bar{d}^{j+\frac{1}{2}}(x_{i})Q_{i}\\&\ge& \frac{\alpha}{\epsilon}Q_{i+1}\bigg(\bar{a}^{j+\frac{1}{2}}(x_{i})-2\alpha\bigg)\\&\geq& \frac{C Q_{i}}{\epsilon+\alpha h_{i}}\ge 0.
\end{eqnarray*}
For large enough $C$, $\eta^{+(j+1)}_{N/2}\leq 0,\eta^{+(0)}_{i}\leq 0$ and $\eta^{+(j+1)}_{N}\leq 0$.
So, by discrete maximum principle, we obtain  \[\eta^{+(j+1)}_{i} \leq 0 \implies |W^{+(j+1)}(x_{i})| \leq C\prod_{n=\frac{N}{2}+1}^{i}\bigg(1+\frac{\alpha h_{n}}{\epsilon}\bigg)^{-1}, \; \frac{N}{2}+1\leq i \leq N-1,~0\le j\le M-1.\]

\end{proof}
The truncation error at the mesh point $(x_{i},t_{j+1})\in\Omega^{N,M}\backslash\{d\}$, when mesh is arbitrary as given by:
\begin{align*}
    |\mathcal{L}_{h}^{N,M}(e(x_{i}))|\le\left\{\begin{array}{ll}
\displaystyle |(\mathcal{L}_{c}^{N,M}-\mathcal{L})y_{i}|\le\epsilon\hbar_{i}\|y_{xxx}\|+\hbar_{i}\|a\|\|y_{xx}\|+\|y_{tt}\|\\
\displaystyle |(\mathcal{L}_{m}^{N,M}-\mathcal{L})y_{i}|\le\epsilon\hbar_{i}\|y_{xxx}\|+C h_{i+1}^2(\|y_{xxx}\|+\|y_{xx}\|)+\|y_{tt}\|
\end{array}
\right.
\end{align*}
On a uniform mesh with step size h, we get
\begin{align*}
    |\mathcal{L}_{h}^{N,M}(e(x_{i}))|\le\left\{\begin{array}{ll}
\displaystyle |(\mathcal{L}_{c}^{N,M}-\mathcal{L})y_{i}|\le\epsilon h^{2}\|y_{xxxx}\|+ h^{2}|a\|\|y_{xxx}\|+\|y_{tt}\|\\
\displaystyle |(\mathcal{L}_{m}^{N,M}-\mathcal{L})y_{i}|\le\epsilon h\|y_{xxx}\|+C h^2(\|y_{xxx}\|+\|y_{xx}\|)+\|y_{tt}\|,
\end{array}
\right.
\end{align*}
where $C_{\|a\|,\|a'\|}$ is a positive constant depends on $\|a\|$ and $\|a'\|$.
\begin{lemma}
\label{regular part}
	The discrete regular component $V^{j+1}(x_{i})$  defined in \eqref{LV} and $ v(x,t)$ is solution of the problem \eqref{v}. So, the error in the regular component satisfies the following estimate:
	$$\|V-v\|_{\Omega^{N-}\cup\Omega^{N+}}\le C(N^{-2}+\Delta t^{2}).$$
\end{lemma}
\begin{proof}
    The truncation error for the regular part of the solution $y$ of the problem (\eqref{oneparaparabolic}), when mesh is uniform $(\tau_{1}=\frac{d}{2})$ in the domain $(x_{i},t_{j+1})\in\Omega^{N-}$ is
\begin{align*}
\lvert \mathcal{L}_{h}^{N,M}(V^{-(j+1)}-v^{-(j+1)})(x_{i})\rvert &=\vert \mathcal{L}_{c}^{N,M}(V^{-(j+1)}-v^{-(j+1)})(x_{i})\rvert \\
&\le \bigg\lvert \epsilon \bigg(\delta^{2}-\frac{d^{2}}{dx^{2}}\bigg)v^{-(j+1)}(x_{i})\bigg\rvert + \lvert a^{j+\frac{1}{2}}(x_{i})\rvert  \bigg\lvert \bigg(D^{0}-\frac{d}{dx}\bigg)v^{-(j+1)}(x_{i})\bigg\rvert\\&+\bigg\lvert \bigg(D_{t}^{-}-\frac{\delta}{\delta t}\bigg)v^{-(j+1)}(x_{i})\bigg\rvert\\
&\le\epsilon h^{2}\|v^{-}_{xxxx}\|+ h^{2}\|a\|\|v^{-}_{xxx}\|+C_{1}\Delta t^{2}\\&\le C(N^{-2}+\Delta t^{2}).\end{align*}
When mesh is non uniform in $(x_{i},t_{j+1})\in\Omega^{N-}$:
\begin{align*}
\lvert \mathcal{L}_{c}^{N,M}(V^{-(j+1)}-v^{-(j+1)})(x_{i})\rvert&\le\epsilon \hbar_{i}\|v^{-}_{xxx}\|+ \hbar_{i}\|a\|\|v^{-}_{xx}\|+C_{1}\Delta t^{2} \le C(N^{-2}+\Delta t^{2}).
\end{align*}
Define the barrier function in $(x_{i},t_{j+1})\in\Omega^{N-}:$
$$\psi^{j+1}(x_{i})=-C(N^{-2}+\Delta t^2)\pm(V^{-(j+1)}-v^{-(j+1)})(x_{i}).$$
For large C,  $\psi^{j+1}(0)\le0, \; \psi^{j+1}(x_{\frac{N}{2}})\le0,\psi^{0}(x_{i})\le0$ and $\mathcal{L}_{h}^{N,M}\psi^{j+1}(x_{i})\ge0$. 
Hence, using the approach given in \cite{FHMRS1}, we get $\psi^{j+1}(x_{i})\le0$ and 
\begin{equation}
\label{V-}
\lvert (V^{-(j+1)}-v^{-(j+1)})(x_i)\rvert_{\Omega^{N-}} \le C(N^{-2}+\Delta t^{2}).
\end{equation}
Similarly, we can estimate the error bounds for right regular part in the domain $(x_{i},t_{j+1})\in\Omega^{N+}$:
\begin{equation}
\label{V+}
\lvert (V^{+(j+1)}-v^{+(j+1)})(x_i)\rvert_{\Omega^{N+}} \le C(N^{-2}+\Delta t^{2}).
\end{equation}
Combining the above results \eqref{V-} and \eqref{V+}, we obtain
$$\|V-v\|_{\Omega^{N-}\cup\Omega^{N+}}\le C(N^{-2}+\Delta t^{2}).$$
\end{proof}
\begin{lemma}
	\label{left layer}
	Let $W^{-}(x_{i},t_{j+1}), w^{-}(x,t)$ are solution of the problem \eqref{LW-} and \eqref{lw} in $\Omega^{-}$ respectively. The left layer component of the truncation error satisfies the following estimate:
 \begin{equation}
     \lvert \|W^{-}-w^{-}\|_{\Omega^{N-}} \le \left\{
	\begin{array}{ll}
	\displaystyle C(N^{-2}+\Delta t^2), & \hbox{ $0\le i\le \frac{N}{4}$, } \vspace{.5mm}\\
\displaystyle	C(N^{-2}(\ln N)^2+\Delta t^2), & \hbox{ $\frac{N}{4}+1\le i \le \frac{N}{2}-1$.}
	\end{array}
	\right.
 \end{equation}
\end{lemma}
\begin{proof}
    The truncation error in case of uniform $(\tau_{1}=\frac{d}{2})$ in the domain $(x_{i},t_{j+1})\in\Omega^{N-}$:
\begin{align*}
    \lvert \mathcal{L}_{c}^{N,M}(W^{-(j+1)}-w^{-(j+1)})(x_{i})\rvert
&\le CN^{-2}(\epsilon\|w_{xxxx}^{-(j+1}\|+\|w_{xxx}^{-(j+1)}\|)+C_{1}\Delta t^{2}\| w^{-}_{tt}\|\\
&\le C(N^{-2}+\Delta t^2).
\end{align*}

In case of non-uniform mesh, we split the argument into two parts.
In first part $(x_{i},t_{j+1})\in (0, d-\tau_{1}]\times(0,T]$ from theorem (\ref{bounds1}), we obtain
	\begin{equation}
	\label{wl}
	\lvert w^{-(j+1)}(x_{i})\rvert \le C\exp^{-(d-x_{i})\alpha_{1}/\epsilon}\le C\exp^{-\tau_{1}\alpha_{1}/\epsilon }\le CN^{-2}.
	\end{equation}
	Also from lemma (\ref{discrete singular}),
	\begin{equation}
 \label{W-}
	\lvert W^{-(j+1)}(x_{i})\rvert \le CR_{i}\le C R_{\frac{N}{4}}=  C\prod_{n=\frac{N}{4}+1}^{\frac{N}{2}}\bigg(1+\frac{\alpha h_{n}}{\epsilon}\bigg)^{-1}\le\bigg(1+\frac{16 \ln N}{N}\bigg)^{-N/8}\le CN^{-2}.
	\end{equation}
 Hence, for all $(x_{i},t_{j+1})\in (0,d-\tau_{1}]\times(0,T]$ by using equation (\ref{wl}) and (\ref{W-}), we have
	$$\lvert (W^{-(j+1)}-w^{-(j+1)})(x_{i})\rvert \le \lvert W^{-(j+1)}(x_{i})\rvert +\lvert w^{-(j+1)}(x_{i})\rvert \le CN^{-2}.$$
The  truncation error for the left layer component  in the inner region $(d-\tau_{1},d)\times(0,T],$ is
 \begin{align*}
\lvert \mathcal{L}_{c}^{N,M}(W^{-(j+1)}-w^{-(j+1)})(x_{i})\rvert
&\le \bigg\lvert \epsilon \bigg(\delta^{2}-\frac{\delta^{2}}{\delta x^{2}}\bigg)w^{-(j+1)}(x_{i})\bigg\rvert + \lvert a^{j+\frac{1}{2}}(x_{i})\rvert  \bigg\lvert \bigg(D^{0}-\frac{\delta}{\delta x}\bigg)w^{-(j+1)}(x_{i})\bigg\rvert\\&+\bigg\lvert \bigg(D_{t}^{-}w^{-(j+1)}(x_{i})-\frac{\delta}{\delta t}w^{-(j+\frac{1}{2})}(x_{i})\bigg)\bigg\rvert\\&\le CN^{-2}(\epsilon\|w_{xxxx}^{-}\|+\|w_{xxx}^{-}\|)+C_{1}\Delta t^2\\
&\le C(N^{-2}(\ln N)^{2}+\Delta t^2).
\end{align*}
We choose the barrier function for the layer component as
	$$\psi^{j+1}(x_{i})=-C_{1}((N^{-1}\ln N)^{2}+\Delta t^2)\pm(W^{-(j+1)}-w^{-(j+1)})(x_{i}),~~(x_{i},t_{j+1})\in(d-\tau_{1},d)\times(0,T].$$
	For sufficiently large $C$, we have $\psi^{j+1}(x_{0})\le0,\psi^{j+1}(x_{N/2})\le0,\psi^{0}(x_{i})\le0,~x_{i}\in(d-\tau_{2},d)$ and $\mathcal{L}_{c}^{N,M}\psi^{j+1}(x_{i})\ge0,~(d-\tau_{2},d)\times(0,T]$.  Hence by discrete maximum principle in \cite{RPS}, $\psi^{j+1}(x_{i})\le0, ~(d-\tau_{2},d)\times(0,T]$.
 So, we can obtain the following bounds:
	$$\lvert (W^{-(j+1)}-w^{-(j+1)})(x_{i})\rvert \le C(N^{-2}(\ln N)^2+\Delta t^2),\,\quad  \forall \; (x_{i},t_{j+1})\in (d-\tau_{2},d)\times(0,T].$$
 So, we can obtain the following bounds in $(0,d)\times(0,T]$:
	\begin{equation}
      \|W^{-}-w^{-}\|_{\Omega^{N-}} \le \left\{
	\begin{array}{ll}
	\displaystyle C(N^{-2}+\Delta t^2), & \hbox{ $0\le i\le \frac{N}{4}$, } \vspace{.5mm}\\
\displaystyle	C(N^{-2}(\ln N)^2+\Delta t^2), & \hbox{ $\frac{N}{4}+1\le i \le \frac{N}{2}-1$.}
	\end{array}
	\right.
 \end{equation}

\end{proof}
\begin{lemma}
	\label{right layer}
	Let $W^{+}(x_{i},t_{j+1}), w^{+}(x,t)$ are solution of the problem \eqref{LW+} and \eqref{lw} in $\Omega^{+}$respectively. The right layer component of the truncation error satisfies the following estimate:
 \begin{equation}
     \|W^{+}-w^{+}\|_{\Omega^{N+}} \le \left\{
	\begin{array}{ll}
\displaystyle	C(N^{-2}(\ln N)^2+\Delta t^2), & \hbox{ $\frac{N}{2}+1\le i \le \frac{3N}{4}-1$,} \vspace{.5mm}\\
\displaystyle C(N^{-2}+\Delta t^2), & \hbox{ $\frac{3N}{4}\le i\le N$. }
	\end{array}
	\right.
 \end{equation}
\end{lemma}
\begin{proof}
  The truncation error when the mesh is uniform $(\tau_{2}=\frac{1-d}{2})$ in $(x_{i},t_{j+1})\in\Omega^{N+}$:
\begin{align*}
    \lvert \mathcal{L}_{c}^{N,M}(W^{+(j+1)}-w^{+(j+1)})(x_{i})\rvert
&\le CN^{-2}(\epsilon\|w_{xxxx}^{-}\|+\|w_{xxx}^{-}\|)+C_{1}\Delta t^2\\
&\le C(N^{-2}+\Delta t^2).
\end{align*}
In case of non-uniform mesh, we split the interval $(d,1)\times(0,T]$ into two parts $(d,d+\tau_{2})\times(0,T]$ and $[d+\tau_{2},1)\times(0,T]$. \\
In the domain $(x_{i},t_{j+1})\in[d+\tau_{2},1)\times(0,T]$:
\begin{equation}
	\label{rightlayerlemma}
	\lvert w^{+(j+1)}(x_{i})\rvert \le C\exp^{-(x_{i}-d)\alpha_{2}/\epsilon}\le C\exp^{- \tau_{2}\alpha_{2}/\epsilon}\le CN^{-2}.
	\end{equation}
	Also from Lemma (\ref{discrete singular}),
	\begin{equation}
 \label{W+}
	    \lvert W^{+(j+1)}(x_{i})\rvert \le C\prod_{n=1}^{\frac{N}{8}}\bigg(1+\frac{\alpha h_{n}}{\epsilon}\bigg)^{-1}\le\bigg(1+\frac{16 \ln N}{N}\bigg)^{-N/8}\le CN^{-2}.
	\end{equation}
 Hence, for all $(x_{i},t_{j+1})\in [d+\tau_{2},1)\times(0,T]$ by using (\ref{rightlayerlemma}) and (\ref{W+}), we have
	$$\lvert (W^{+(j+1)}-w^{+(j+1)})(x_{i})\rvert \le \lvert W^{+(j+1)}(x_{i})\rvert +\lvert w^{+(j+1)}(x_{i})\rvert \le CN^{-2}.$$
 The  truncation error for the left layer component  in the inner region $(d, d+\tau_{2})\times(0,T],$ is
 \begin{align*}
\lvert \mathcal{L}_{c}^{N,M}(W^{+(j+1)}-w^{+(j+1)})(x_{i})\rvert
&\le \bigg\lvert \epsilon \bigg(\delta^{2}-\frac{d^{2}}{\delta x^{2}}\bigg)w^{+(j+1)}(x_{i})\bigg\rvert + \lvert a^{j+\frac{1}{2}}(x_{i})\rvert  \bigg\lvert \bigg(D^{0}-\frac{\delta}{\delta x}\bigg)w^{+(j+1)}(x_{i})\bigg\rvert\\&+\bigg\lvert \bigg(D_{t}^{-}w^{-(j+1)}(x_{i})-\frac{\delta}{\delta t}w^{+(j+\frac{1}{2})}(x_{i})\bigg)\bigg\rvert\\&\le CN^{-2}(\epsilon\|w_{xxxx}^{+}\|+\|w_{xxx}^{+}\|)+C_{1}\Delta t^{2}\\
&\le C(N^{-2}(\ln N)^{2}+\Delta t^2).
\end{align*}
We choose the barrier function for the layer component as
	$$\psi^{j+1}(x_{i})=-C_{1}((N^{-1}\ln N)^{2}+\Delta t^2)\pm(W^{+(j+1)}-w^{+(j+1)})(x_{i}),~~(x_{i},t_{j+1})\in(d, d+\tau_{2})\times(0,T].$$
	For sufficiently large $C$, we have $\psi^{j+1}(d)\le0,\psi^{j+1}(x_{\frac{3N}{4}})\le0,\psi^{0}(x_{i})\le0,x_{i}\in (d, d+\tau_{2})$ and $\mathcal{L}_{h}^{N,M}\psi^{j+1}(x_{i})\ge0,(x_{i},t_{j+1})\in(d, d+\tau_{2})\times(0,T]$.  Hence by discrete maximum principle in \cite{RPS}, $\psi^{j+1}(x_{i})\le0$.
 So, we can obtain the following bounds:
	$$\lvert (W^{+(j+1)}-w^{+(j+1)})(x_{i})\rvert \le C(N^{-2}(\ln N)^2+\Delta t^2),\,\quad  \forall \; (x_{i},t_{j+1})\in(d, d+\tau_{2})\times(0,T].$$
 Hence, the desired bound is:
	\begin{equation}
    \|W^{+}-w^{+}\|_{\Omega^{N+}} \le \left\{
	\begin{array}{ll}
\displaystyle	C(N^{-2}(\ln N)^2+\Delta t^2), & \hbox{ $\frac{N}{2}+1\le i \le \frac{3N}{4}-1$,} \vspace{.5mm}\\
\displaystyle C(N^{-2}+\Delta t^2), & \hbox{ $\frac{3N}{4}\le i\le N$. }
	\end{array}
	\right.
 \end{equation}
\end{proof}

\begin{lemma}
		The error $e(x_{\frac{N}{2}},t_{j+1})$ estimated at the point of discontinuity $(x_{\frac{N}{2}},t_{j+1})=(d,t_{j+1}),0\le j\le M-1$ satisfies the following estimates:
		$$\displaystyle \lvert \mathcal{L}_{h}^{N,M}(Y-y)(x_{\frac{N}{2}},t_{j+1})\rvert \le C\bigg(\frac{ \tau^2}{\epsilon^{3}N^2}+\Delta t^2\bigg)$$
\end{lemma}
\begin{proof}
    The  truncation error at the point $x_{\frac{N}{2}}=d$, where $\tau=\tau_{1}=\tau_{2}$:
\begin{align*}
	\bigg\lvert \mathcal{L}_{h}^{N,M} Y^{j+1}(x_{N/2})-\frac{h g^{j+1}_{N/2-1}}{2\epsilon-h a^{j+\frac{1}{2}}_{N/2-1}}-\frac{h g^{j+1}_{N/2+1}}{2\epsilon+h a^{j+\frac{1}{2}}_{N/2+1}}\bigg\rvert &\le\lvert\mathcal{L}_{t}^{N,M} Y^{j+1}(x_{N/2})\rvert+\frac{h}{2\epsilon-h a^{j+\frac{1}{2}}_{N/2-1}}\\&\lvert\mathcal{L}_{c}^{N,M} y(x_{N/2-1},t_{j+1})-g^{j+1}_{N/2-1}\rvert +\frac{h}{2\epsilon+h a^{j+\frac{1}{2}}_{N/2+1}}\\&\lvert\mathcal{L}_{c}^{N,M} y(x_{N/2+1},t_{j+1})-g^{j+1}_{N/2+1}\rvert\\&\le Ch^{2}\|y_{xxx}\|_{\Omega^{-}\cup\Omega^{+}} \\&\le C\bigg(\frac{\tau^2}{\epsilon^{3}N^2}+\Delta t^2\bigg)
	\end{align*}
\end{proof}
\begin{theorem}
	\label{main}
    Let $y(x,t)$ and $ Y(x,t)$ be the solutions of problem \eqref{oneparaparabolic} and \eqref{discrete problem} respectively, then
	$$\|Y-y\|_{\bar{\Omega}^{N,M}}\le
	C(N^{-2}\ln N^3+\Delta t^{2})$$
	where C is a constant independent of $\epsilon, $ and discretization parameter $N,M$. 
\end{theorem}
\begin{proof}
From Lemma (\ref{regular part}), Lemma (\ref{left layer}), and Lemma (\ref{right layer}), we have that
	 	$$\displaystyle \|Y-y\|_{\Omega^{N-}\cup\Omega^{N+}} \le C(N^2\ln N^{-2}+\Delta t^2)$$
  To find error at the point of discontinuity $x_{\frac{N}{2}}=d$, we have considered the discrete barrier function $\phi(x_{i},t_{j+1})=\psi(x_{i},t_{j+1})\pm e(x_{i},t_{j+1})$ defined in the interval $(d-\tau_{1},d+\tau_{2})\times(0,T]$, where
	 $$\psi(x_{i},t_{j+1})=-C(N^{-2}(\ln N)^{3}+\Delta t^2)-\frac{C \tau^2}{\epsilon^{3}N^2}\left\{
	 \begin{array}{ll}
	 	\displaystyle x_{i}-(d+\tau), & \hbox{ $(x_{i},t_{j+1})\in \Omega^{N,M}\cap(d-\tau, d)\times(0,T]$, }\\
	 	d-\tau-x_{i}, & \hbox{ $(x_{i},t_{j+1})\in\Omega^{N,M}\cap(d, d+\tau)\times(0,T].$}
	 \end{array}
	 \right. $$
	 We have $\phi(d-\tau,t_{j+1})$, $\phi(d+\tau,t_{j+1}),\phi(x_{i},t_{0}),~x_{i}\in(d-\tau,d+\tau) $ are non positive and
	 $\mathcal{L}_{h}^{N,M}\phi(x_{i},t_{j+1})\ge 0,~~(x_{i},t_{j+1})\in( d-\tau, d+\tau)\times(0,T],~~ \mathcal{L}_{t}^{N,M}\phi(x_{N/2},t_{j+1}) \ge 0.$\\
  Hence, by applying discrete maximum principle we get $\phi(x_{i},t_{j+1})\le0.$\\
	Therefore, for $(x_{i},t_{j+1})\in (d-\tau,d+\tau)\times(0,T]:$
	\begin{equation}
	\label{Yd}
	\lvert (Y-y)(x_{i},t_{j+1})\rvert \le C_{1}(N^{-2}(\ln N)^3+\Delta t^2)+C_{2}\frac{h^{2}}{\epsilon^{3}}\tau+\Delta t^2)\le C(N^{-2}(\ln N)^3+\Delta t^2).
	\end{equation}
\end{proof}
\section{Numerical examples} 
To demonstrate the effectiveness of the proposed hybrid difference approach, we conducted numerical experiments on two test problems with discontinuous convection coefficients and source terms. Here, for both the examples $\Omega= (0,1)\times (0,1], d=0.5,$ and $ N=M$.
\begin{Example}\label{ex-a}
	$$(\epsilon y_{xx}+ ay_{x}-by-y_{t})(x,t)=f(x,t),~~ ~~(x,t)\in((0,.5)\cup(0.5,1))\times(0,1],$$
	$$	y(0,t)=y(1,t)=y(x,0)=0,$$\nonumber\
	with \\
	$$a(x,t)=\left\{
	\begin{array}{ll}
	\displaystyle -(1+x(1-x)), & \hbox{ $0\le x\le 0.5,t\in(0,1]$, }\\
	1+x(1-x), & \hbox{ $0.5<x\le1,t\in(0,1]$,}
	\end{array}
	\right.$$
	$$
	f(x,t)=\left\{
	\begin{array}{ll}
	\displaystyle -2(1+x^2)t, & \hbox{ $0\le x\le 0.5,t\in(0,1]$, }\\
	2(1+x^2)t, & \hbox{ $0.5<x\le1,t\in(0,1]$,}
	\end{array}
	\right.$$\\
	and $b(x,t)=1+exp(x)$.
\end{Example}
\begin{Example}\label{ex-b}
	$$(\epsilon y_{xx}+ ay_{x}-by-y_{t})(x,t)=f(x,t),~~ ~~(x,t)\in((0,.5)\cup(0.5,1))\times(0,1],$$
	$$	y(0,t)=y(1,t)=y(x,0)=0,$$\nonumber\
	with \\
	$$a(x,t)=\left\{
	\begin{array}{ll}
	\displaystyle -1, & \hbox{ $0\le x\le 0.5,t\in(0,1]$, }\\
	1, & \hbox{ $0.5<x\le1,t\in(0,1]$,}
	\end{array}
	\right.$$
	$$
	f(x,t)=\left\{
	\begin{array}{ll}
	\displaystyle -2xt, & \hbox{ $0\le x\le 0.5,t\in(0,1]$, }\\
	2(1-x)t, & \hbox{ $0.5<x\le1,t\in(0,1]$,}
	\end{array}
	\right.$$\\
	and $b(x,t)=0$.
\end{Example}
As the exact solutions to these problems are unknown, we used the double mesh method \cite{DV,D} to evaluate the accuracy of the numerical approximations obtained. The double mesh difference is defined by
$$E_{\epsilon}^{N,M}= \max\limits_{j}\bigg(\max\limits_{i}|Y_{2i,2j}^{2N,2M}-Y_{i,j}^{N,M}|\bigg)$$
where $Y_{i}^{N,M}$  is the solution on the mesh $\bar{\Omega}^{N,M}$ and $Y_{2i}^{2N,2M}$ is the solutions on the mesh  $\bar{\Omega}^{2N,2M}$ respectively where the fine mesh $\bar{\Omega}^{N,M}$ and  $\bar{\Omega}^{2N,2M}$ are with $N$ and $2N$ mesh-intervals in the spatial direction and $2M$ mesh-intervals in the
t-direction.
The order of convergence is given by
$$R_{\epsilon}^{N,M}=\log_{2}\bigg(\frac{E_{\epsilon}^{N,M}}{E_{\epsilon,}^{2N,2M}}\bigg).$$
\begin{table}
\small
\begin{center}
\caption{Maximum point-wise error $E_{\epsilon}^{N,M}$ and approximate orders of convergence $R_{\epsilon}^{N,M}$ for Example \ref{ex-a}.}\label{tab1}%
\begin{tabular}{cccccc}
\hline
\multirow{2}{*}{$\epsilon$ }& \multicolumn{5}{c}{Number of mesh points N}\\

   & 64 &128 & 256&512&1024\\
\hline
		$2^{-8}$ & 1.12e-01&	4.07e-02&	1.39e-02&	4.53e-03&	1.42e-03 \\ 
		Order &1.4645&	1.5483&	1.6182&	1.6698& \\

		$2^{-10}$ &1.14e-01&	4.09e-02&	1.40e-02&	4.55e-03&	1.43e-03\\ 
		Order & 1.4822&	1.5409&	1.6261&	1.6699&\\
		$2^{-12}$ &1.15e-01&	4.09e-02&	1.41e-02&	0.004.56e-03&	1.43e-03 \\ 
		Order & 1.4853&	1.5402&1.6265&	1.6691& \\
  
		$2^{-14}$ & 1.15e-01&	4.10e-02&	1.41e-02&	4.57e-03&	1.44e-03\\ 
		
		Order & 1.4863&	1.5407&	1.6250&	1.6709& \\
		
		$2^{-16}$ &1.15e-01&	4.10e-02&	1.41e-02&	4.57e-03&	1.45e-03\\
		Order & 1.4837&	1.5407&	1.6250&	1.6709& \\
		$2^{-18}$ &1.15e-01&	4.10e-02&	1.41e-02&	4.57e-03&	1.45e-03\\
		Order & 1.4837&	1.5407&	1.6250&	1.6709& \\
		$2^{-20}$ &1.15e-01&	4.10e-02&	1.41e-02&	4.57e-03&	1.45e-03\\
		Order & 1.4837&	1.5407&	1.6250&	1.6709& \\
   \vdots &\vdots&\vdots&\vdots&\vdots&\vdots\\
  $2^{-40}$ &1.15e-01&	4.10e-02&	1.41e-02&	4.57e-03&	1.45e-03\\
		Order & 1.4837&	1.5407&	1.6250&	1.6709& \\
\hline
\end{tabular}
\end{center}
\end{table}

\begin{table}
\small
\begin{center}
\caption{Maximum point-wise error $E_{\epsilon}^{N,M}$ and approximate orders of convergence $R_{\epsilon}^{N,M}$ for Example \ref{ex-a}.}\label{tab2}%
\begin{tabular}{cccccc}
\hline
\multirow{2}{*}{$\epsilon$ }& \multicolumn{5}{c}{Number of mesh points N}\\
   & 64 &128 & 256&512&1024\\
\hline
		$2^{-8}$ & 5.04e-02&	1.84e-02&	6.44e-03&	2.13e-03&6.73e-04\\ 
		Order &1.4533&	1.5148&	1.5933&	1.6642& \\

		$2^{-10}$ &4.90e-02&	1.79e-02&	6.28e-03&	2.07e-03&	6.50e-04\\ 
		Order & 1.4521&	1.5132&	1.5990&	1.6729&\\
		$2^{-12}$ &4.87e-02&	1.76e-02&	6.17e-03&	2.03e-03&6.35e-04 \\ 
		Order &1.4699	1.5115&	1.6039&	1.6755& \\
  
		$2^{-14}$ &4.86e-02&1.75e-02&6.15e-03&	2.02e-03&	6.34e-04\\ 
		
		Order & 1.4699&	1.5138&	1.6032&	1.6735& \\
		
		$2^{-16}$ &4.86e-02&1.76e-02&6.15e-3& 2.02e-03&	6.34e-04\\
		Order & 1.4676&	1.5105&	1.6032&	1.6735& \\
		$2^{-18}$ &4.86e-02&1.76e-02&	6.15e-03&	2.02e-03&	6.34e-04\\
		Order & 1.4676&	1.5106&	1.6032&	1.6735& \\
		$2^{-20}$ &4.86e-02&1.76e-02&	6.15e-03&	2.02e-03&	6.34e-04\\
		Order & 1.4676&	1.5106&	1.6032&	1.6735& \\
   \vdots &\vdots&\vdots&\vdots&\vdots&\vdots\\
  $2^{-40}$ &4.86e-02&1.76e-02&	6.15e-03&	2.02e-03&	6.34e-04\\
		Order & 1.4676&	1.5106&	1.6032&	1.6735& \\
\hline
\end{tabular}
\end{center}
\end{table}
 \begin{figure}[h]
\centering
\includegraphics[width=0.5\textwidth]{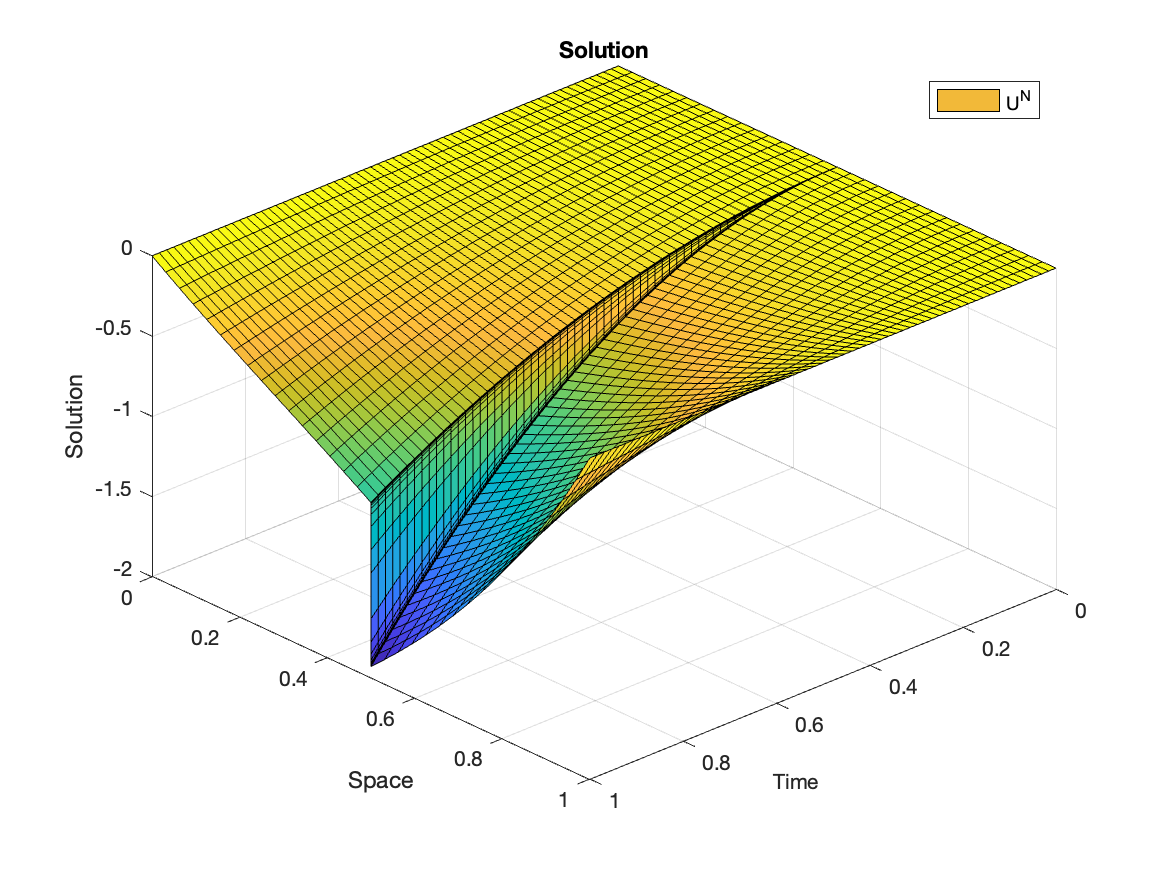}
\caption{Numerical solution for $\epsilon=2^{-16} $when $N=64$ for Example \ref{ex-a}.}	
\label{fig1}
\end{figure}

\begin{figure}[h]
\centering
\includegraphics[width=0.5\textwidth]{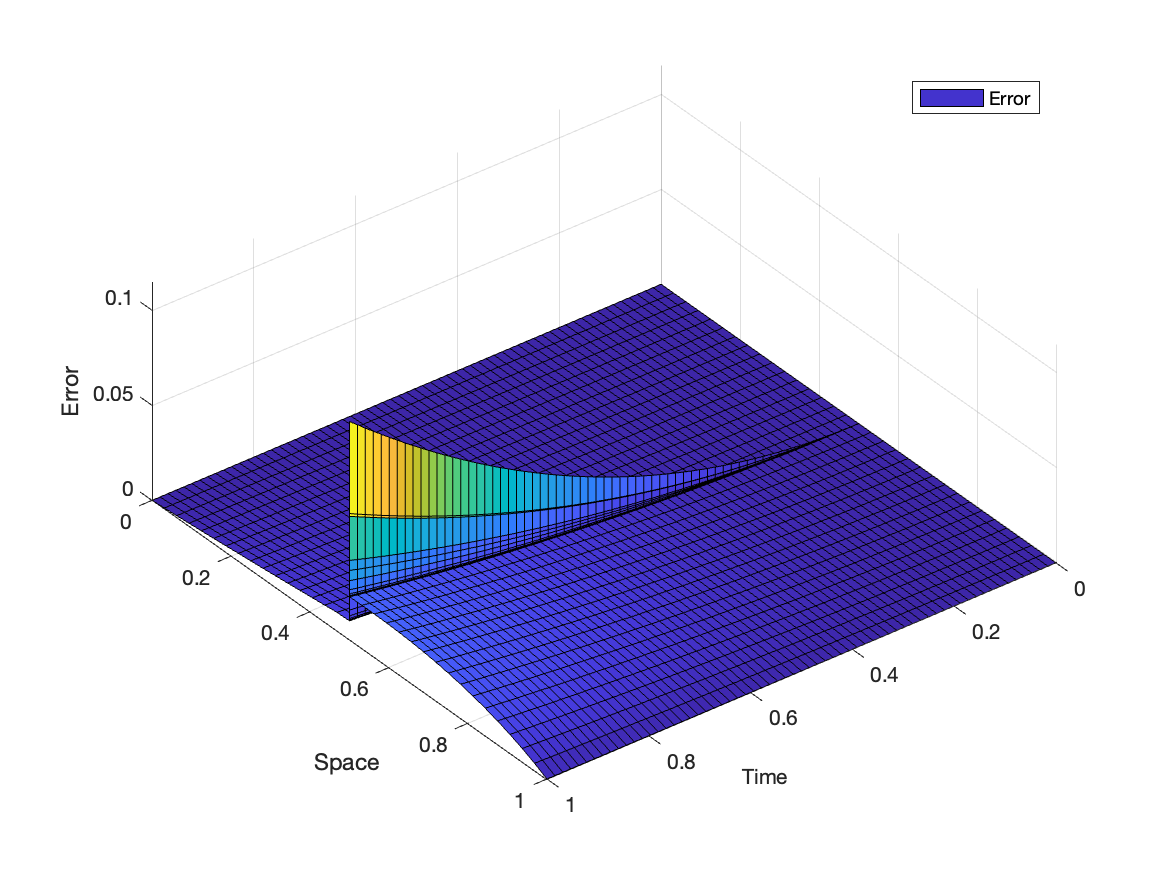}
\caption{Error for $\epsilon=2^{-16}$ when $N=64$ for Example \ref{ex-a}.}	
\label{fig2}
\end{figure}
\begin{figure}[h]
\centering
\includegraphics[width=0.5\textwidth]{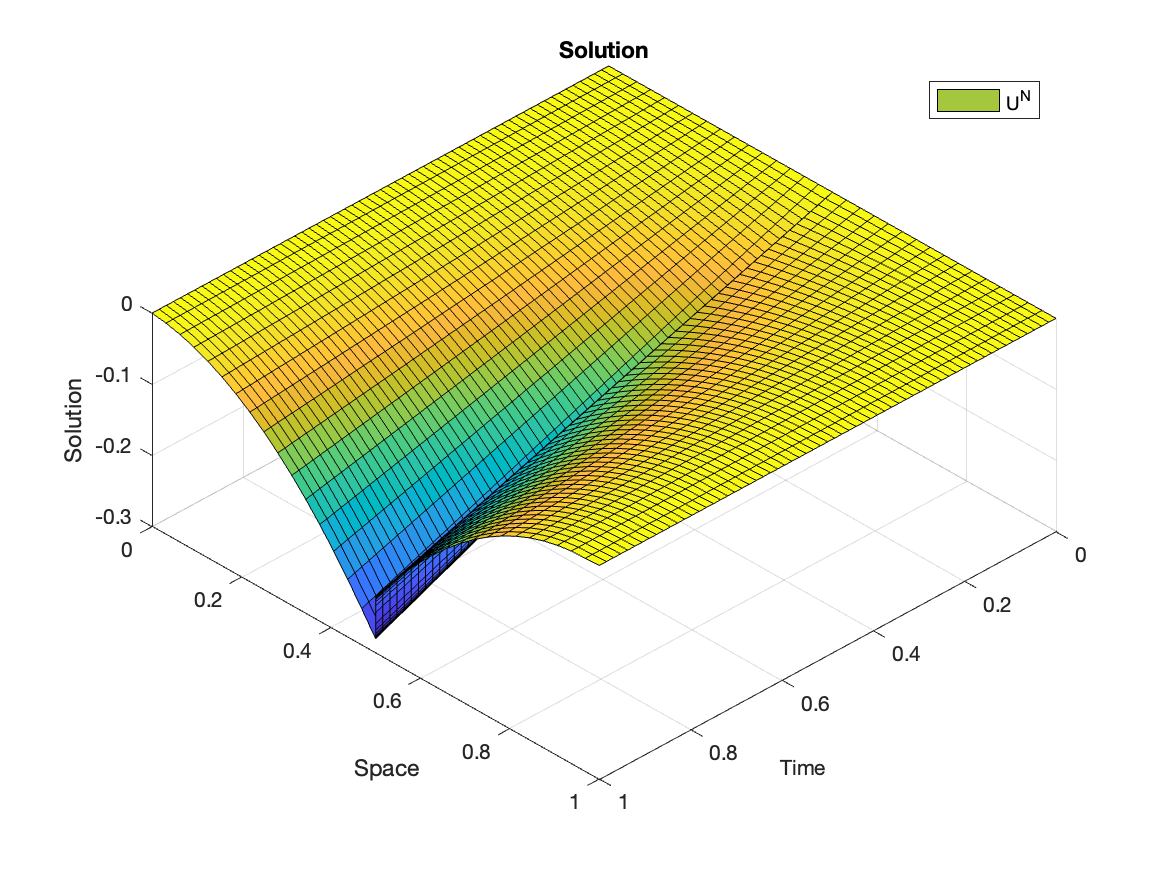}
\caption{Numerical solution for $\epsilon=2^{-22}$ when $N=64$ for Example \ref{ex-b}.}	
\label{fig3}
\end{figure}
\begin{figure}[h]
\centering
\includegraphics[width=0.5\textwidth]{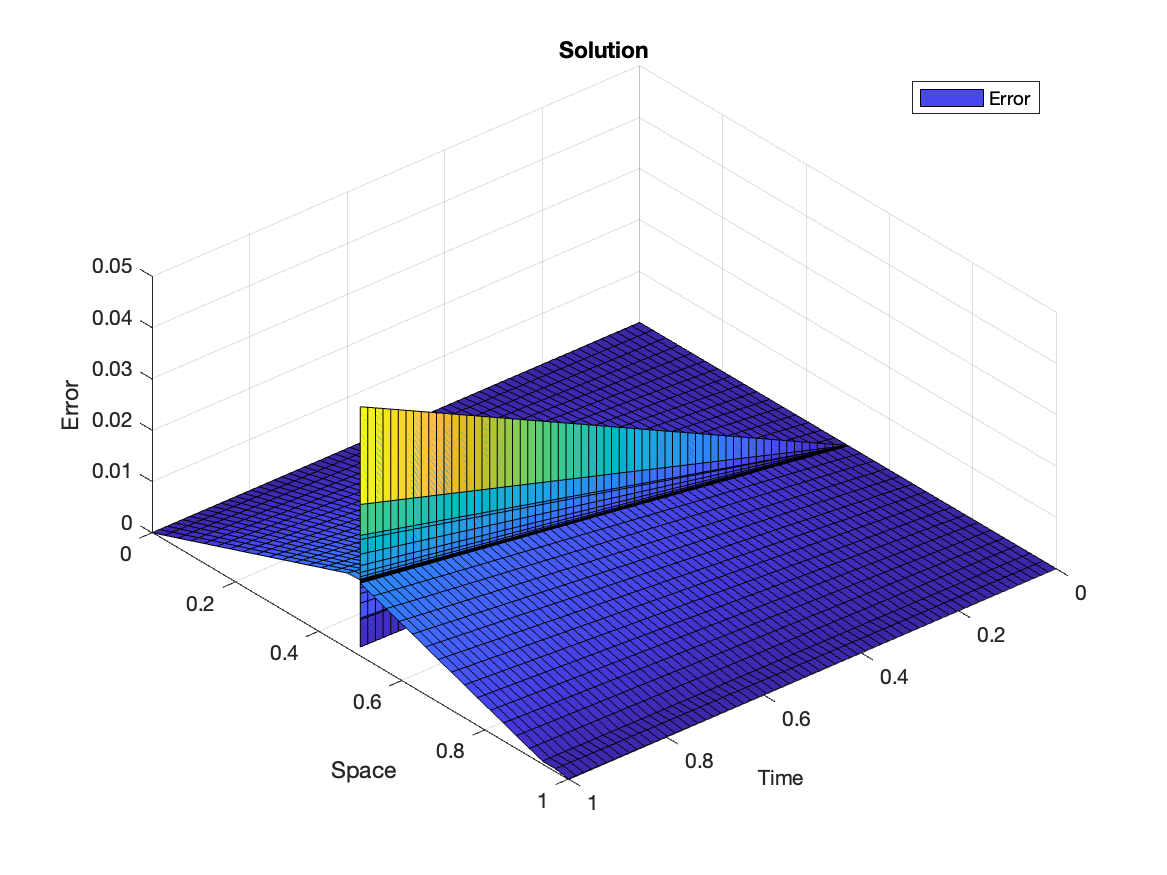}
\caption{Error for $\epsilon=2^{-22}$ when $N=64$ for Example \ref{ex-b}.}	
\label{fig4}
\end{figure}
In table (\ref{tab1}) and table (\ref{tab2}), we have presented the maximum pointwise  error and the corresponding order of convergence for different values of $\epsilon$ for example (\ref{ex-a}) and (\ref{ex-b}), as $\epsilon$ increases, layers becomes more steep and pronounced. Also $E_{\epsilon}^{N,M}$ decreases monotonically as N increases and the order of convergence is increases up to almost second order. 
We have also plotted the solution and error graph in figure (\ref{fig1}-\ref{fig2}) and figure (\ref{fig3}-\ref{fig4}) for both the example (\ref{ex-a}) and (\ref{ex-b}) for the values of $\epsilon=2^{-16}$ and $\epsilon=2^{-22}$ when $N=64$. We can observe from the figure that the maximum error is obtained at the point of discontinuity.
\section{Conclusion}
In this study, we investigated a singularly perturbed one-parameter parabolic equation with a discontinuous convection coefficient and a source term. The solution exhibits internal layers at the points of discontinuity, caused by the discontinuities in the convection coefficient and source term. For the temporal discretization, we applied the Crank-Nicolson method on a uniform mesh, achieving second-order accuracy in time. In the spatial direction, we employed a hybrid difference scheme that combines the midpoint method with the central difference method, using a specially designed Shishkin mesh and a five-point formula at the discontinuity. This approach ensures nearly second-order accuracy in space whish is better than any other results. Theoretical findings were supported by numerical results, confirming the effectiveness and accuracy of the proposed scheme.

\end{sloppypar}
\end{document}